\documentclass[12pt]{amsart}
\usepackage{amsmath,amsfonts,amsthm,amscd}
\newtheorem{theorem}{Theorem}[section]
\newtheorem{claim}[theorem]{Claim}
\newtheorem{proposition}[theorem]{Proposition}
\newtheorem{lemma}[theorem]{Lemma}
\newtheorem{definition}[theorem]{Definition}

\theoremstyle{remark}
\newtheorem*{remark}{Remark}
\numberwithin{equation}{section}
\begin{document}
\bibliographystyle{amsalpha}
\title{Moduli spaces of critical Riemannian metrics in dimension four} 
\author{Gang Tian}
\address{Gang Tian\\
 Department of Mathematics \\ MIT \\ Cambridge, MA 02139}
\curraddr{Department of Mathematics, Princeton University, 
Princeton, NJ 08544}
\email{tian@math.mit.edu}
\thanks{The research of the first author was 
partially supported by NSF Grant DMS-0302744.}
\author{Jeff Viaclovsky}
\address{Jeff Viaclovsky, Department of Mathematics, MIT, Cambridge, MA 02139}
\email{jeffv@math.mit.edu}
\thanks{The research of the second author was partially 
supported by NSF Grant DMS-0202477.}
\date{December 16,2003}
\begin{abstract}
We obtain a compactness result for various classes of 
Riemannian metrics in dimension four; in particular our
method applies to anti-self-dual metrics, K\"ahler metrics with 
constant scalar curvature, and metrics with harmonic curvature. 
With certain geometric non-collapsing assumptions, the moduli space can be 
compactified by adding metrics with orbifold-like singularities. Similar results 
were obtained for Einstein metrics in \cite{Anderson}, \cite{BKN}, 
\cite{Tian}, but our analysis differs substantially from the 
Einstein case in that we do not assume any pointwise Ricci curvature bound.
\end{abstract}
\maketitle
\section{Introduction}
 Critical points of the total scalar curvature functional
(restricted to the space of unit volume metrics)
\begin{align}
\label{functional2}
\mathcal{R} : g \mapsto \int_M R_g \ dV_g, 
\end{align}
are exactly the Einstein metrics, and the 
structure of the moduli space of Einstein metrics
has been extensively studied (\cite{Anderson}, \cite{BKN}, 
\cite{Nakajima2}, \cite{Tian}). In particular, with certain 
geometric assumptions on non-collapsing, this moduli space can 
be compactified by adding Einstein metrics with 
orbifold singularities. 

   The motivation for this paper is to prove 
a similar compactness theorem for various 
classes of metrics in dimension four, where 
one does not assume a pointwise bound
on the Ricci curvature.  We will consider the 
following cases:\\
\begin{tabular}{ll}
\\
a. & half-conformally flat metrics constant scalar curvature,\\

b. & metrics with harmonic curvature,\\

c. & K\"ahler metrics with constant scalar curvature.\\ 
\end{tabular}
\\
\\
Half conformally flat metrics are also known 
as self-dual or anti-self-dual if $W^- = 0$
or $W^+ = 0$, respectively. These metrics  
are, in a certain sense, analogous to anti-self-dual 
connections in Yang-Mills theory (see \cite{FreedUhlenbeck},
\cite{DonaldsonKronheimer}). 
The local structure of the moduli space of anti-self-dual 
metrics, by examining the linearization of the
anti-self-dual equations, has been studied, for 
example, in \cite{AHS}, \cite{Itoh} and \cite{Kotschick}.   
There has been a considerable amount of research on the 
existence of anti-self-dual metrics on compact manifolds.
In the paper \cite{Poon}, Poon constructed  
a one-parameter family of anti-self-dual conformal classes on
$ \overline{\mathbb{CP}}^2 \# \overline{\mathbb{CP}}^2$.
LeBrun \cite{LeBrun2} produced explicit 
examples on $n \overline{\mathbb{CP}}^2 $
for all $n \geq 2$. We also mention the work of 
\cite{Floer} and \cite{DonaldsonFriedman} for other methods 
and examples. See also \cite{LeBrun5} 
for a nice survey and further references. A very important 
contribution is Taubes' stable existence theorem for anti-self-dual 
metrics: for any compact, oriented, smooth four-manifold M, the manifold
$M \# n \overline{\mathbb{CP}}^2$
carries an anti-self-dual metric for some $n$ (see \cite{Taubes}). 
This shows that anti-self-dual metrics exists in 
abundance, so one would like to understand the 
moduli space.

 In \cite{Bo4}, it was proved that a compact four-dimensional 
Riemannian manifold with harmonic 
curvature and non-zero signature must be Einstein. 
Therefore (b) is larger than the class of
Einstein metrics only in the case of zero signature.
In particular, we have locally conformally flat metrics 
with constant scalar curvature, which have been studied 
in \cite{SchoenYau1}, \cite{SchoenYau2},
\cite{Schoen3}.
For more background about cases (a)-(c) above, 
see \cite{Besse}. We also note that case (c) 
is an important class of extremal K\"ahler metrics
\cite{Calabi1}, \cite{Calabi2}.

  In the sequel, when we say {\em{critical metric}} we 
will mean any of (a)-(c) above. 
For $M$ compact, we define the Sobolev constant $C_S$  
as the best constant $C_S$ so that 
for all $f \in C^{0,1}(M)$ we have
\begin{align}
\label{mainSob}
\Vert f \Vert_{L^{4}} \leq C_S   
\Vert \nabla f \Vert_{L^2} + Vol^{-1/4} \Vert f \Vert_{L^2},
\end{align}
where $Vol$ is the volume. Note that (\ref{mainSob}) is scale-invariant. 
For $M$ non-compact, $C_S$ is defined to be the best constant so that 
\begin{align}
\label{mainSob2}
\Vert f \Vert_{L^{4}} \leq C_S   
\Vert \nabla f \Vert_{L^2},
\end{align}
for all $f \in C^{0,1}(M)$ with compact support. 

  We define a {\em{Riemannian 
orbifold}} $(M,g)$ to be a topological space which is a 
smooth manifold with a smooth Riemannian metric 
away from finitely many singular points.  
At a singular point $p$, $M$ is locally diffeomorphic 
to a cone $\mathcal{C}$ on 
$S^3 / \Gamma$, where $\Gamma \subset SO(4)$ 
is a finite subgroup acting freely on $S^3$.
Furthermore, at such a singular point, the metric is locally the 
quotient of a smooth $\Gamma$-invariant metric on $B^4$ under 
the orbifold group $\Gamma$.
We note that the notions of smooth orbifold,
orbifold diffeomorphism, and orbifold Riemannian 
metric are well-defined, see 
\cite{Satake1}, \cite{Satake2}, \cite{Thurston}, \cite{Borz1},
\cite{TianYau} for more background.
A Riemannian orbifold $(M,g)$ 
is a {\em{K\"ahler orbifold}} if
$g$ is K\"ahler, all of the orbifold groups $\Gamma$
are in $U(2)$, and at each singular point, the metric is 
locally the quotient of a smooth K\"ahler metric on a ball 
in $\mathbf{C}^2$ under the orbifold group.

Consider the disjoint union 
\begin{align}
\label{odecomp}
\tilde{M} = \coprod_{i=1}^{N} M_{i},
\end{align} 
where each $M_{i}$ is a Riemannian
orbifold. Then a {\em{Riemannian multi-fold}} $M$ 
is a connected space obtained from $\tilde{M}$ by finitely many identifications 
of points. Note that points from $M_i$ and $M_j$, $i \neq j$
can be identified, as well as several points from the 
same $M_i$. For example, take $M_1$ and $M_2$ to be smooth manifolds, 
and identify $p_1 \in M_1$ with $p_2 \in M_2$. 
Another example would be to take just one smooth manifold $M_1$, 
and identify $p_1 \in M_1$ with $p_2 \in M_1$.
The {\em{singular set}} of $M$ is the set of  
points where $M$ is not a smooth manifold -- this will come 
from the nontrivial orbifold singular points of each $M_i$, 
as well as new singular points from the identifications.
These latter points look like multiple cone points, thus the 
terminology {\em{multi-fold}}.  
If there is more than one orbifold in (\ref{odecomp})
($N > 1$), some $M_{i}$ is compact, and has 
only one point which gets identified to the 
other orbifolds $M_j$, $i \neq j$ to form $M$, then 
we say $M$ {\em{splits off}} the compact orbifold $M_i$. 
If there is only one cone at a singular point 
$p$, then $p$ is called {\em{irreducible}}. 

 A smooth Riemannian manifold $(M,g)$ 
is called an asymptotically locally 
Euclidean (ALE) end of order $\tau$ 
if there exists a finite subgroup $\Gamma \subset SO(4)$ 
acting freely on $\mathbf{R}^4 \setminus B(0,R)$ and a 
$C^{\infty}$ diffeomorphism 
$\Psi : M \rightarrow ( \mathbf{R}^4 \setminus B(0,R)) / \Gamma$ 
such that under this identification, 
\begin{align}
g_{ij} &= \delta_{ij} + O( r^{-\tau}),\\
\ \partial^{|k|} g_{ij} &= O(r^{-\tau - k }),
\end{align}
for any partial derivative of order $k$ as
$r \rightarrow \infty$. 
We say an end 
is ALE of order $0$ if we can find a coordinate system 
as above with $g_{ij} = \delta_{ij} + o(1)$, 
and $\partial^{|k|} g_{ij} = o(r^{- k })$
as $r \rightarrow \infty$. 
A complete, noncompact Riemannian multi-fold $(M,g)$ is called ALE 
if $M$ can be written as the disjoint union of a compact 
set and finitely many ALE ends. 

We say that a sequence of Riemannian manifolds 
$(M_j, g_j)$ {\em{converges}} 
to the Riemannian multi-fold $(M_{\infty}, g_{\infty})$ 
if the following is satisfied. 
For $\epsilon > 0$, consider $M_{\infty, \epsilon}
= M_{\infty} \setminus S_{\epsilon}$, where $S_{\epsilon}$
is the $\epsilon$-neighborhood of $S$, and $S$ is a finite 
set of points containing all of the singular points of
$M_{\infty}$.
Then there exist domains $\Omega_j(\epsilon) \subset M_j$, 
and diffeomorphisms $\Phi_{j, \epsilon}: M_{\infty, \epsilon}
\rightarrow \Omega_j(\epsilon)$, such that 
$\Phi_{j, \epsilon}^{\ast} g_j$ converges to 
$g_{\infty}$ in $C^{\infty}$ as $j \rightarrow \infty$,
on compact subsets of $M_{\infty, \epsilon}$. 
Furthermore, there exist constants 
$\delta, N$ depending upon $\epsilon$,
such that 
\begin{align}
\max \{ Vol( M_j \setminus \Omega_j(\epsilon)), 
diam(  M_j \setminus \Omega_j(\epsilon))\} < \delta,
\end{align}
for $j > N$ and 
$\delta \rightarrow 0$ as $\epsilon \rightarrow 0$,
where $Vol$ and $diam$ denote the volume and diameter 
with respect to the metric $g_j$, respectively. 
A sequence of pointed Riemannian manifolds 
$(M_j, g_j, p_j)$ converges 
to the pointed Riemannian multi-fold $(M_{\infty}, g_{\infty}, p_{\infty})$ 
if for all $R > 0$, $B(p_j,R) $ converges to $B(p_{\infty}, R)$ 
as above as pointed spaces. 

 We state our main convergence theorem:  
\begin{theorem} 
\label{orbthm}
Let $(M_i, g_i, p_i)$ be a sequence of critical metrics $g_i$ on 
smooth, complete, pointed four-dimensional manifolds $M_i$ satisfying
\begin{align} 
\label{bound1}
&C_S \leq C_1, \\
\label{l2curv}
&\int_{M_i} |Rm_{g_i}|^2 dV_{g_i} \leq \Lambda,
\end{align}
\begin{align}
\label{diamvol}
&  Vol(g_i) > \lambda > 0,\\
\label{betti1}
& b_1(M_i) < B_1, 
\end{align}
where $C_1, \Lambda, \lambda$  are constants,
and $b_1(M_i)$ denotes the first Betti number. 
Then there exists a subsequence $\{ j \} \subset \{ i \}$, 
a pointed, connected, critical 
Riemannian multi-fold $(M_{\infty}, g_{\infty}, p_{\infty})$, 
and a finite singular set $S \subset M_{\infty}$ such that 
\\
(1) $(M_j, g_j,p_j)$ converges to $(M_{\infty}, g_{\infty}, p_{\infty})$.
\\
(2) The limit space $(M_{\infty}, g_{\infty}, p_{\infty})$ 
does not split off any compact orbifold.
\\
(3) If $M_{\infty}$ is non-compact,  then $(M_{\infty}, g_{\infty}, 
p_{\infty})$ is ALE of order $\tau$ for any $\tau < 2$.
\\
(4)  If $b_1(M_i) = 0$, then $(M_{\infty}, g_{\infty}, p_{\infty})$ is 
a Riemannian orbifold.
\\
(5) In the  K\"ahler case (c), $(M_{\infty}, g_{\infty}, p_{\infty})$
is a K\"ahler orbifold. 
\end{theorem}
\begin{remark}
We note that the definition of convergence given here 
implies, in particular, Gromov-Hausdorff convergence. 
Moreover, we will show in Section~\ref{secconv} that the convergence is 
even stronger, in the sense that suitable rescalings of the metrics 
near the singular points converge to ALE multi-folds. 
 The singular set $S$ is the singular set of 
convergence, it necessarily contains the multi-fold 
singular set of $M_{\infty}$, but 
it is possible for some points in $S$ to be smooth 
points of $M_{\infty}$. This is in contrast to the
Einstein case, where the Bishop-Gromov volume 
comparison theorem implies that 
convergence is smooth at any smooth point in $M_{\infty}$. 
\end{remark}
\begin{remark}
A positive lower 
bound on the Yamabe invariant $Y(M_i, [g_i])$ will imply the Sobolev
constant bound, and in certain geometric situations, this bound 
will be automatically satisfied.
The bound in (\ref{l2curv}) will also follow automatically 
in certain geometric situations. 
We will discuss these points in Section \ref{sobvconst} below.
Furthermore, the main elements of our proof only require a 
{\em{local}}  Sobolev constant bound, see Theorem \ref{orbthmsob} below. 
\end{remark}

 In conjunction with Theorem \ref{orbthm}, we have the volume comparison theorem:
\begin{theorem} 
\label{orbthm2}
Let $(M, g)$ be a critical metric on a 
smooth, complete four-dimensional manifold $M$ satisfying
\begin{align} 
&C_S \leq C_1,\\
&\int_{M} |Rm_{g}|^2 dV_g \leq \Lambda,\\
& b_1(M) < B_1, 
\end{align}
where $C_1, \Lambda,$ and $B_1$ are constants.
Then there exists a constant 
$V_1$, depending only upon $C_1, \Lambda, $ and $B_1$, 
such that $Vol (B(p,r)) \leq V_1 \cdot r^4$, 
for all $p \in M$ and $r > 0$.  
\end{theorem}
Finally, we restate Theorem \ref{orbthm} in the compact case:
\begin{theorem}
\label{orbthm3}
Let $(M_i, g_i)$ be a sequence of critical metrics on 
smooth, closed four-dimensional manifolds $M_i$ satisfying
\begin{align} 
&C_S \leq C_1,\\
&\int_{M_i} |Rm_{g_i}|^2 dV_{g_i} \leq \Lambda,\\
& Vol(M_i,g_i)=1,\\
& b_1(M_i) < B_1, 
\end{align}
where $C_1, \Lambda, B_1$ are constants.
Then the conclusion of Theorem \ref{orbthm} 
holds. That is, the limit space $(M_{\infty}, g_{\infty})$ is a compact,
connected, critical Riemannian multi-fold which does not split off
any compact orbifold. In the K\"ahler case (c), or if $b_1(M_i) =0$, 
then $M_{\infty}$ is an orbifold.  
\end{theorem} 
\begin{remark}
All of our results hold in the more general Bach-flat case
(see Section \ref{criticalequation}), with the exception 
that at a singular point, we can only show the metric 
$g$ is locally the quotient by the orbifold group of a 
$C^0$ metric on a standard ball, smooth away from the origin,
and in the ALE case, the metric is ALE of order $0$. 
\end{remark}

\subsection{Acknowledgements}

 The authors benefitted on several occasions from 
conversations with Denis Auroux, Tom Branson, Toby Colding, John Lott,
Tom Mrowka, Jie Qing, and Joao Santos. In particular, we are indebted to 
Tom Branson for assistance with the Kato 
inequality in Section \ref{removal}. 
We would also like to thank the referee for valuable 
comments which helped to improve the exposition 
of the paper. 

 This paper is the second part of \cite{TV}.
We mention the recent preprint \cite{Andersonc}, which 
deals with similar issues. The work 
here was done independently, but the first author would like to acknowledge 
discussions with Michael Anderson on this subject prior to 1995. 
\section{Critical metrics}
\label{criticalequation}
In this section, we discuss the systems of equations  
satisfied in cases (a), (b), and (c), and 
justify the terminology {\em{critical metric}}. 

\subsection{Half-conformally flat metrics and metrics with harmonic curvature}

 These systems were discussed in \cite{TV}, so 
we just briefly review them here. 
 
 The Euler-Lagrange equations of the functional
\begin{align}
\label{functional3}
\mathcal{W} : g \mapsto \int_M |W_g|^2 dV_g, 
\end{align}
in dimension four, are
\begin{align}
\label{bacheq2}
B_{ij} = \nabla^k \nabla^l W_{ikjl} +
\frac{1}{2}R^{kl}W_{ikjl} = 0.
\end{align}
where $W_{ijkl}$ and $R_{kl}$ are the components 
of the Weyl  and Ricci tensors, respectively (see \cite{Besse}, 
\cite{Derdzinski}). 
Since the Bach tensor arises in the Euler-Lagrange 
equations of a Riemannian functional, it is symmetric,
and since the functional (\ref{functional3}) is 
conformally invariant, it follows that the Bach-flat 
equation (\ref{bacheq2}) is conformally invariant.
The Bach tensor arises as the Yang-Mills equation 
for a twistor connection \cite{Merkulov},
see also \cite{BastonMason}, \cite{LeBrun3} for 
other occurrences of the Bach tensor. 

We note that (see \cite{Calderbankweaklyasd})
\begin{align}
B_{ij} = 2 \nabla^k \nabla^l W^+_{ikjl} +
R^{kl}W^+_{ikjl} = 2 \nabla^k \nabla^l W^-_{ikjl} +
R^{kl}W^-_{ikjl},
\end{align}
so that both self-dual and anti-self-dual metrics are 
Bach-flat. 

Using the Bianchi identities,
a computation shows that we may rewrite the 
Bach-flat equation (in dimension four) as
\begin{align}
(\Delta Ric)_{ij} = 2 ( R_{il}g_{jk} - R_{ikjl} -  W_{ikjl})
( R_{kl} - (R/6) g_{kl}). 
\end{align}
Introducing a convenient shorthand, we write this as 
\begin{align}
\label{shortcurveqn}
\Delta Ric = Rm* Ric.
\end{align}

The condition for harmonic curvature is 
that 
\begin{align}
\label{harmcond}
\delta Rm = -R_{ijkl;i} = 0. 
\end{align}
This condition was studied in \cite{Bo4}, 
\cite{Derd1}, \cite{Besse}, and is 
the Riemannian analogue of a Yang-Mills 
connection. An equivalent condition for harmonic curvature
that $\delta W = 0$ and $ R = constant$. 
In particular, locally conformally flat 
metrics with constant scalar curvature have 
harmonic curvature. A computation shows (\ref{shortcurveqn}) 
is also satisfied, but in this case we moreover 
have an equation on the full curvature tensor. We compute
\begin{align*}
(\Delta Rm)_{ijkl} &= R_{ijkl;m;m}\\
&= ( - R_{ijlm;k} - R_{ijmk;l})_{;m}\\
&= -R_{ijlm;mk}  - R_{ijmk;ml} + Q(Rm)_{ijkl} = Q(Rm)_{ijkl}.
\end{align*}
where $Q(Rm)$ denotes a quadratic expression 
in the curvature tensor. 
In the shorthand, we write this as 
\begin{align}
\label{harmcurveqn}
\Delta Rm = Rm* Rm.
\end{align}

\subsection{K\"ahler metrics with constant scalar 
curvature} 
We assume that $(M,g)$ is a K\"ahler manifold with 
K\"ahler metric $g$. In \cite{Calabi1} it was shown that if 
$dR$ is a holomorphic vector field, then $g$ is critical 
for the $L^2$ norm of the scalar curvature, 
restricted to a K\"ahler class. In particular K\"ahler
and constant scalar curvature implies extremal. 

 The bisectional curvature tensor is given by 
\begin{align*}
R_{i \overline{j} k \overline{l}} = - \frac{ \partial^2 g_{i \overline{j}}}
{ \partial z_k \partial \overline{z}_l}
+ g^{s \overline{t}} \frac{\partial g_{s \overline{j}}}{\partial \overline{z}_l}
 \frac{\partial g_{i \overline{t}}}{\partial z_k},
\end{align*}
in local coordinates $(z_1, \dots, z_n)$ of $M$. 
Contracting with the inverse of $\{ g_{i \overline{j}} \}$,
we obtain for the Ricci and scalar curvature
\begin{align*}
R_{i \overline{j}} &= - \frac{ \partial^2}{ \partial z_i \partial \overline{z}_j}
( \log \det ( g_{k \overline{l}})),\\
R &= -2 \Delta_g \log \det (  g_{k \overline{l}})
= - 2 \frac{ \partial^2}{ \partial z_k \partial \overline{z}_k}
( \log \det ( g_{i \overline{j}})).
\end{align*}
In particular these imply the following Bianchi identities
\begin{align*}
R_{ i \overline{j} k \overline{l};m} &= R_{ i \overline{j} m \overline{l} ; k}\\
R_{i \overline{j}; \overline{k}} &= R_{i \overline{k} ; \overline{j}}\\
R_{i \overline{j} k \overline{l}} &= R_{k \overline{j} i \overline{l}}.
\end{align*}
It follows then that in local unitary frames 
\begin{align*}
\Delta_g( Ric)_{i \overline{j}} &= R_{i \overline{j}; k ; \overline{k}}\\
&= R_{k \overline{j}; i ; \overline{k}} \\
&=  R_{k \overline{j}; \overline{k} ; i} + R_{i \overline{s}}R_{s \overline{j}}
- R_{k \overline{s}} R_{i \overline{j} s \overline{k}}\\
&= \frac{1}{2} R_{;i ; \overline{j}} + R_{i \overline{s}}R_{s \overline{j}}
- R_{k \overline{s}} R_{i \overline{j} s \overline{k}}.
\end{align*}
Therefore if the scalar curvature is constant, we have 
\begin{align}
\begin{split}
\Delta_g( Ric)_{i \overline{j}} &=  R_{i \overline{s}}R_{s \overline{j}}
- R_{k \overline{s}} R_{i \overline{j} s \overline{k}}\\
& = (R_{ i \overline{s}} g_{k \overline{j}} 
 -  R_{k \overline{s}} R_{i \overline{j} s \overline{k}})
R_{k \overline{s}},
\end{split}
\end{align}
so in this case, we again have the equation 
\begin{align}
\Delta Ric = Rm* Ric.
\end{align}
\section{Geometric bounds}
\label{sobvconst}
In this section, we will discuss some special cases 
for which various assumption in Theorem 
\ref{orbthm} will be automatically satisfied.
We recall that the Yamabe invariant in dimension four is 
defined by 
\begin{align*}
Y(M,[g_0]) = \underset{g \in [g_0]}{ \mbox{inf}} Vol(g)^{-1/2}
\int_M R_g dV_g.
\end{align*}
We define the Sobolev constant as the best constant 
such that for all $\phi \in C^{0,1}_c(M)$,
\begin{align*}
\Vert \phi \Vert_{L^{4}}
 \leq C_S  \Vert \nabla \phi \Vert_{L^2} + Vol^{-1/4}\Vert 
\phi \Vert_{L^2} . 
\end{align*}
\begin{proposition}
If $g$ is a Yamabe minimizer, and $Y(M,[g]) > 0$, then 
$C_S(M,g) \leq \sqrt{6} Y(M,[g])^{-1/2}$.
\end{proposition}
\begin{proof}
From the definition of the Yamabe invariant, for any $u \in L^2_1(M)$, 
\begin{align*}
\int_M \Big(6 |\nabla_g u|_g^2 + R_g u^2 \Big) dV_g 
\geq Y(M) \left\{ \int_M u^4 dV_g \right\}^{1/2},
\end{align*}
where we use $g$ as the background metric. Since $g$ has 
constant scalar curvature, this implies
\begin{align*}
\frac{6}{Y(M)} \int_M |\nabla u|^2 + \frac{R_g Vol(g)^{1/2}}{Y(M)} 
Vol(g)^{-1/2}\int_M u^2 
\geq \left\{ \int_M u^4 \right\}^{1/2}.
\end{align*}
Since $g$ is Yamabe, we have $R_g Vol(g)^{1/2}=  Y(M)$, so 
we obtain
\begin{align*}
\Vert u \Vert_{L^{4}} \leq  \sqrt{6} Y(M)^{-1/2}
\Vert \nabla u \Vert_{L^2} + Vol(g)^{-1/4} \Vert f \Vert_{L^2}.
\end{align*}
\end{proof}
In dimension four, the Gauss-Bonnet and signature formulas are (see \cite{Besse})\begin{align}
\label{Euler}
8 \pi^2 \chi(M) &= \frac{1}{6}\int_M R^2 - \frac{1}{2} \int_M |Ric|^2 + 
\int_M |W|^2,\\
\label{Signature}
12 \pi^2 \tau(M) &=  \int_M |W^+|^2 - \int_M |W^-|^2.
\end{align}
 In the anti-self-dual case, $W^+ \equiv 0$, we have 
\begin{align}
\label{asd1}
8 \pi^2 \chi(M) &= \frac{1}{6}\int_M R^2 - \frac{1}{2}\int_M |Ric|^2 + \int_M |W^-|^2,\\
\label{asd2}
12 \pi^2 \tau(M) &=  - \int_M |W^-|^2.
\end{align}
Since the anti-self-dual equation is conformally invariant, 
we can make a conformal change to a Yamabe minimizer (
see \cite{Aubin}, \cite{Schoen}, \cite{LeeandParker}), and add 
these equations together to obtain 
\begin{align}
\label{added}
8 \pi^2 \chi(M) +12 \pi^2 \tau(M)  &= \frac{R^2}{6}Vol(M) - \frac{1}{2}\int_M |Ric|^2.
\end{align}
If $R > 0$, and $2 \chi(M) + 3 \tau(M) > 0$ then we obtain the estimate 
\begin{align} 
\label{added2}
Y(M, [g]) = R Vol(M)^{1/2} \geq 2 \sqrt{6} \pi ( 2 \chi + 3 \tau) > c
>0. 
\end{align}
\begin{proposition}
Let $(M,g)$ be Yamabe with $R >0$, and  
anti-self-dual. Then $\Vert Rm_g \Vert_{L^2} < C$, where $C$ 
depends only on $\chi(M), \tau(M)$. 
Furthermore, if  $2 \chi(M) + 3 \tau(M) > 0$, then 
the Sobolev constant is uniformly bounded from above,  
\begin{align}
2 \pi (C_S)^2 \leq \sqrt{6} \left( 2 \chi(M) + 3 \tau(M)\right)^{-1}.
\end{align} 
\end{proposition}
\begin{proof}
For the first statement, (\ref{added}) gives 
a bound on $\Vert Ric \Vert_{L^2}$, since 
$Y(M,[g]) \leq Y(S^n)$ (see \cite{LeeandParker}). 
Equation (\ref{asd2})  gives a bound on $\Vert W \Vert_{L^2}$.
The second statement follows from (\ref{added2}).
\end{proof}

We next consider the K\"ahler case. 
Let $c_1(M)$ denote the first Chern class of $M$. 
It is known that for complex surfaces,
\begin{align*}
c_1^2(M) = 2 \chi(M) + 3 \tau(M),
\end{align*} 
and therefore on a complex surface, 
\begin{align}
\label{cherncond}
Q'(M, [g]) \equiv c_1^2(M)  
- \frac{1}{3}\frac{(c_1(M) \cdot \omega_g)^2}{\omega_g^2}
\end{align}
is a conformal invariant. It follows that  
when $Q'(M, [g]) >0$,
\begin{align}
Y(M, [g]) \geq 3 \pi^2
\sqrt{ Q'(M,[g]) }.
\end{align}
This implies 
\begin{proposition}
For $(M,g)$ K\"ahler satisfying
\begin{align} 
\label{sobv}
3 c_1^2(M) > (c_1(M) \cdot [\omega_g]) ^2,
\end{align}
the Sobolev constant is uniformly bounded from above.  
\end{proposition}

\subsection{On the Sobolev inequality}
All of the results in this paper are still valid 
if the weaker Sobolev inequality is assumed
\begin{align}
\label{weakSob}
\Vert \phi \Vert_{L^{4}}
 \leq C_S  \left( \Vert \nabla \phi \Vert_{L^2} + Vol^{-1/4}\Vert 
\phi \Vert_{L^2} \right), 
\end{align}
with the exceptions being 
in (2) in Theorem \ref{orbthm}, $M_{\infty}$ might 
split off a compact orbifold, and in (4) of Theorem \ref{orbthm}, 
even if $b_1(M_i) = 0$, the limit may be reducible, 
see the proof of Proposition \ref{finitely1} below. 
Furthermore,  as remarked in the introduction, the main elements of 
our proof only require a {\em{local}} Sobolev constant bound, 
see Theorem \ref{orbthmsob} below.

 If we have a conformal class with positive Yamabe invariant, 
we have shown above that the Sobolev constant of the 
Yamabe minimizer is bounded. However, if we instead
choose a non-minimizing constant scalar 
curvature metric, we will have a Sobolev 
inequality of type (\ref{weakSob}). 

\section{Local regularity}
\label{localregularity}

 In all the above cases, the equation take the 
form 
\begin{align}
\label{generaleqn0}
(\Delta Ric)_{ij} = A_{ikjl}R_{kl},
\end{align}
where $A_{ikjl}$ is some linear expression 
in the curvature tensor. 
Using a convenient shorthand, we write this as 
\begin{align}
\label{generaleqn}
\Delta Ric = Rm * Ric.
\end{align}
Using the Bianchi identities, any Riemannian 
metric satisfies 
\begin{align}
\label{curveqn}
\Delta Rm = L(\nabla^2 Ric) + Rm* Rm,
\end{align}
where $L(\nabla^2 Ric)$ denotes a linear expression 
in second derivatives of the Ricci tensor. 

 Even though second derivatives of the 
Ricci occur in (\ref{curveqn}), overall the principal
symbol of the system (\ref{generaleqn}) and (\ref{curveqn})
in triangular form. The equations (\ref{generaleqn}) and (\ref{curveqn}), 
when viewed as an elliptic system, together with the bound on the 
Sobolev constant, yield the following $\epsilon$-regularity 
theorem:
\begin{theorem}(\cite[Theorem 3.1]{TV})
\label{higherlocalregthm}
Assume that (\ref{generaleqn}) is satisfied, 
choose $r < \mbox{diam}(M)/2$, and let $B(p,r)$ be a geodesic 
ball around the point $p$, and $k \geq 0$. Then there exist  
constants $\epsilon_0, C_k$  (depending upon $C_S$) so that if 
\begin{align*}
\Vert Rm \Vert_{L^2(B(p,r))} = 
\left\{ \int_{B(p,r)} |Rm|^2 dV_g \right\}^{1/2} \leq \epsilon_0,
\end{align*}
then 
\begin{align*}
\underset{B(p, r/2)}{sup}| \nabla^k Rm| \leq
\frac{C_k}{r^{2+k}} \left\{ \int_{B(p,r)} |Rm|^2 dV_g \right\}^{1/2}
\leq \frac{C_k \epsilon_0}{r^{2+k}}. 
\end{align*}
\end{theorem}
\begin{remark} 
We state the following slight variation 
of the above. Let $C_S(r)$ denote the Sobolev constant 
for compactly supported functions in $B(p,r)$, that is, 
\begin{align}
\Vert f \Vert_{L^{4}(B(p,r))}
 \leq C_S(r) \Vert \nabla f \Vert_{L^2(B(p,r))}, 
\end{align}
for all $f \in C^{0,1}_c(B(p,r))$.
Then there exists a universal constant $\epsilon_0$ such that if 
\begin{align}
\left\{ C_S(r)^4 \cdot   
\int_{B(p,r)} |Rm|^2 dV_g \right\}^{1/2} \leq \epsilon_0,
\end{align}
then 
\begin{align*}
\sup_{B(p,r/2)} |Rm| \leq \frac{C}{r^2}   
\left\{ C_S(r)^4  \cdot \int_{B(p,r)} |Rm|^2 dV_g \right\}^{1/2}
\leq \frac{C \epsilon_0}{r^2}
\end{align*}
where $C$ is a universal constant, the proof being the same as 
that of Theorem \ref{higherlocalregthm}. It is also interesting
to bound $C_S(r)$ in terms of the volume of $B(p,r)$. For the 
manifolds being considered in this paper, it may be possible that 
$C_S(r) \cdot Vol(B(p,r))^{1/4} < Cr$, for some uniform 
constant $C$. 
\end{remark}
Theorem \ref{higherlocalregthm} may be applied to 
noncompact orbifolds to give a rate of curvature decay at
infinity. Assume that $(M,g)$ is a complete, noncompact 
orbifold with finitely many singular points, 
with a critical metric, bounded Sobolev constant
(for functions with compact support), and finite $L^2$ norm of 
curvature. Fix a basepoint $p$, and let $r(x) = d(p,x)$. 
Given $\epsilon < \epsilon_0$ from Theorem \ref{higherlocalregthm}, 
there exists an $R$ large so that
there are no singular points on $D(R/2)$ and 
\begin{align*}
\int_{D(R)} |Rm|^2 dV_g < \epsilon < \epsilon_0,
\end{align*}
where $D(R) = M \setminus B(R)$. 
Choose any $x \in M$ with $d(x,p) = r(x) > 2R$, 
then $B(x,r) \subset D(R)$. From 
Theorem \ref{higherlocalregthm}, we have 
\begin{align*}
\underset{B(x, r/2)}{sup}| \nabla^k Rm| \leq
\frac{C_k}{r^{2+k}} \left\{ \int_{B(x,r)} |Rm|^2 dV_g \right\}^{1/2}
\leq \frac{C_k \epsilon}{r^{2+k}}, 
\end{align*}
which implies
\begin{align*}
| \nabla^k Rm|(x)  \leq \frac{C_k \epsilon}{r^{2+k}}.
\end{align*}
As we take $R$ larger, we may choose $\epsilon$ smaller, 
and we see that $M$ has better-than-quadratic curvature decay, along 
with decay of derivatives of curvature. 
\section{Volume Growth}
\label{euclideanvolumegrowth}
One of the crucial aspects of this problem 
is to obtain control on volume growth of metric balls 
from above. We let $(M,d)$ be a length space with 
distance function $d$, and basepoint $p \in M$. 
\begin{definition}
We say a component $A_0(r_1,r_2)$ of an annulus 
$A(r_1, r_2) = \{ q \in M \ | \ r_1 < d(p, q) < r_2 \}$
is {\em{bad}} if $S(r_1) \cap \overline{A_0(r_1,r_2)}$ has 
more than $1$ component, where $S(r_1)$ is the 
sphere of radius $r_1$ centered at $p$. 
\end{definition}
 As we remarked after \cite[Theorem 4.1]{TV}, the proof 
of our volume growth theorem requires only that there
are finitely many disjoint bad annuli, therefore 
we have
\begin{theorem}(\cite[Theorem 4.1]{TV})
\label{bigthm}
Let $(M,g)$ be a complete, non-compact, four-dimensional 
Riemannian orbifold (with finitely many singular
points) with base point $p$.
Assume that there exists a constant $C_1 > 0$ so that
\begin{align}
\label{cond4}
Vol(B(q,s)) \geq C_1 s^4,
\end{align}
for any $q \in M$, and all $s \geq 0$.
Assume furthermore that as $r \rightarrow \infty$,
\begin{align}
\label{decay1}
\underset{S(r)}{sup} \ |Rm_g| &= o(r^{-2}),
\end{align}
where $S(r)$ denotes the sphere of radius $r$ centered
at $p$. If $(M,g)$ contains only finitely many 
disjoint bad annuli, then $(M,g)$ has finitely many ends, 
and there exists a constant $C_2$ so that 
\begin{align}
\label{vga}
Vol(B(p,r)) \leq C_2 r^4,
\end{align}
Furthermore, each end is ALE of order $0$. 
\end{theorem}
\begin{proof}
Since there are no orbifold singular points outside of a compact 
set, the proof of \cite[Theorem 4.1]{TV} is valid in this case.
To see this, from \cite[Proposition 15]{Borz1} any minimizing geodesic segment cannot
pass through the singular set unless it begins or ends 
on the singular set, and the set $ M \setminus S$ is 
geodesically convex. Therefore, all of the standard 
tools from Riemannian geometry used in the 
proof of \cite[Theorem 4.1]{TV} apply in this setting. 
\end{proof}
By taking instead sequences of dyadic annuli 
$A(s^{-j-1}, s^{-j}), 1 < s$, around a singular point, 
the proof of \cite[Theorem 4.1]{TV} can also be applied directly to 
components of isolated singularities:
\begin{theorem}
\label{bigthm2}
Let $(X,d,x)$ be a complete, locally compact length 
space, with basepoint $x$. Let $B(x,1) \setminus \{x\}$
be a $C^{\infty}$ connected four-dimensional manifold with 
a metric $g$ of class (a), (b), or (c) satisfying
\begin{align}
& \int_{B(x,1)} |Rm_g|^2 dV_g < \infty,\\
&  \Vert u \Vert_{L^4( B(x,1) \setminus \{x\})} 
\leq C_s  \Vert \nabla u \Vert_{L^2( B(x,1) \setminus \{x\})}
, \  u \in C^{0,1}(B(x,1) \setminus \{x\}) \\
& b_1(X) < \infty,
\end{align}
where $C_s, V_1$ are positive constants. 
Then there exists a constant $C_1 > 0$ 
so that $Vol (B(x,r)) \leq C_1 r^4.$ 
The basepoint $x$ is an orbifold point, and 
and the metric $g$ extends to $B(x,1)$
as a $C^0$-orbifold metric. 
That is, for some small $\delta >0$, the
universal cover of $B(x, \delta) \setminus \{x\}$ 
is diffeomorphic to a punctured ball $B^4 \setminus \{0\}$ 
in $\mathbf{R}^4$, and the lift of $g$, after diffeomorphism, 
extends to a $C^0$ metric $\tilde{g}$ on $B^4$, 
which is smooth away from the origin. 
\end{theorem}
\begin{remark}
This is valid for components of 
$B(x, \delta) \setminus \{x\}$, 
we will prove below that there 
are finitely many components for the limit 
space arising in Theorem \ref{orbthm}.
To show $x$ is a $C^0$-orbifold point, one uses 
a tangent cone analysis as in \cite[Theorem 4.1]{TV}.
Furthermore, in Theorem \ref{remsing} below, we will show $g$ is 
a {\em{smooth}} orbifold metric. 
\end{remark}
\section{Asymptotic curvature decay and removal of singularities 
with bounded energy}
\label{removal}
We first discuss curvature decay results from \cite[Section 6]{TV},
and using the same technique, we prove a singularity removal result.
\begin{theorem}
\label{decayale}
Let $(M,g)$ be a complete, noncompact four-dimensional irreducible 
Riemannian orbifold with $g$ of class (a), (b), or (c), and 
finitely many singular points. Assume that
\begin{align}
\int_M |Rm_g|^2 dV_g < \infty, \ C_S < \infty,
\mbox{ and }  b_1(M) < \infty.
\end{align}
Then $(M,g)$ has finitely many ends, and 
each end is ALE of order $\tau$ for any $\tau < 2$.
\end{theorem}
\begin{remark} 
In case $(M,g)$ is a manifold, from 
\cite[Theorem 1]{Carron}, we have a bound on the 
number of ends depending only upon the Sobolev constant
and the $L^2$-norm of curvature (moreover, 
all of the the $L^2$-Betti numbers are bounded).
In the K\"ahler case, an argument as in 
\cite{Li-Tam} shows that there can be at 
most 1 non-parabolic end, we remark that 
the analysis there is valid for irreducible 
orbifolds with finitely many singular points. 
Since any ALE end is non-parabolic, this implies there
only one end. The argument in \cite[Theorem 4.1]{Li-Tam}
is roughly,  to construct a nonconstant bounded harmonic function
with finite Dirichlet integral if there is more than 
$1$ non-parabolic end. This function must then be pluriharmonic, 
and under the curvature decay conditions, it must therefore be constant. 
\end{remark}
\begin{proof}
Theorem \ref{decayale} was proved in \cite[Theorem 1.3]{TV}, 
the proof there is also valid for orbifolds.
We briefly outline the details. 
\begin{lemma}
If $(M,g)$ satisfies (a), (b), or (c), then  
\begin{align}
\label{rickato}
\big| \nabla |E| \big|^2 \leq \frac{2}{3} | \nabla E | ^2,
\end{align}
at any point where $|E| \neq 0$, where 
$E$ denotes the traceless Ricci tensor.  
\end{lemma}
This is due to Tom Branson, the proof follows 
from his general theory of Kato constants 
developed in \cite{Branson}, see
\cite[Lemma 5.1]{TV} for the details of this case, 
the proof being valid also in all 
cases (a), (b), and (c). We remark that the same constant follows 
from the methods in \cite{CGH}. The case considered in Lemma 
\ref{rickato} is exactly the case $r=s=2$ in the last line of the table  
on \cite[page 253]{CGH}, giving immediately the best constant $2/3$.

Using this improved Kato constant, we now have the equation
\begin{align} 
\label{impriceqn}
\Delta |E|^{1/2} \geq - C |E|^{1/2} |Rm|. 
\end{align}
Using a Moser iteration argument from \cite{BKN}, 
and since the scalar curvature is constant, this allows 
one to improve the Ricci curvature decay
to $|Ric| = O(r^{-2 - \delta})$ for any $\delta <2$, 
where $r(x) = d(p,x)$ is the distance to a basepoint $p$.
Next, using a Yang-Mills argument (inspired by the proof of Uhlenbeck for 
Yang-Mills connections \cite{Uhlenbeck2}, also \cite[Section 4]{Tian}) 
the following was proved in \cite[Lemma 6.5]{TV}
\begin{lemma}
\label{bali}
Let $D(r) = M \setminus B(p,r)$. 
For $\delta < 2$, and $r$ sufficiently large, we have 
\begin{align}
\underset{D(2r)}{sup} \Vert Rm_g \Vert_g
\leq \frac{C}{r^{2 + \delta}}.
\end{align}
\end{lemma}
The result then follows by the construction of
coordinates at infinity in \cite{BKN}.
\end{proof}

 Next we discuss a removable singularity result,
this is an analogue of \cite[Theorem 5.1]{BKN}, \cite[Lemma 4.5]{Tian}.
This theorem is crucial to obtain smoothness of
the limit orbifold. 
\begin{theorem}
\label{remsing}
Let $(X,d,x)$ be a complete, locally compact length 
space, with basepoint $x$. Let $B(x,1) \setminus \{x\}$
be connected $C^{\infty}$ four-dimensional manifold with 
a metric $g$ of class (a), (b), or (c) satisfying
\begin{align}
& \int_{B(x,1)} |Rm_g|^2 dV_g < \infty,\\
&  \Vert u \Vert_{L^4( B(x,1) \setminus \{x\})} 
\leq C_s  \Vert \nabla u \Vert_{L^2( B(x,1) \setminus \{x\})}
, \  u \in C^{0,1}(B(x,1) \setminus \{x\}) \\
& Vol (B(x,r)) \leq V_1 r^4, \  r > 0,
\end{align}
where $C_s, V_1$ are positive constants. 
Then the metric $g$ extends to $B(x,1)$
as a smooth orbifold metric. 
That is, for some small $\delta >0$, 
universal cover of $B(x, \delta) \setminus \{x\}$ 
is diffeomorphic to a punctured ball $B^4 \setminus \{0\}$ 
in $\mathbf{R}^4$, and the lift of $g$, after diffeomorphism, 
extends to a smooth critical metric $\tilde{g}$
on $B^4$.
\end{theorem}
\begin{proof}
The argument  in \cite[Lemma 6.5]{TV} for ALE spaces examined 
the behavior at infinity, 
we now imitate the argument using balls around a singular point. 
From Theorem \ref{bigthm2} above, we know the singularity 
is orbifold of order $0$, and the tangent 
cone at a singularity is a cone on a spherical space 
form $S^3 / \Gamma$,  We lift by the action of the 
orbifold group to obtain a critical metric in 
$B(0,\delta) \setminus \{ 0 \}$ with bounded energy, bounded Sobolev 
constant, and $Vol(B(0,s)) < Cs^4$. 
From the Kato inequalities in cases (a), (b), and (c), 
we obtain the estimate 
$|Ric| = O(r^{-2 + \delta}),$ where $r$ now denotes 
distance to the origin, for any $\delta < 2$.  
The argument from \cite[Lemma 6.5]{TV} shows 
that $|Rm| = O(r^{-2 + \delta})$. 
As in \cite[Lemma 4.4]{Tian}, we can then find a self-diffeomorphism  
$\psi$ of $B(0, \delta)$ so that $\nabla (\psi^{*} g) = O(r^{-1+\delta}),$ and 
$ \psi^{*} g = O(r^{ \delta})$. Choosing $\delta$ close to $2$,
the metric $\psi^{*}g$ has a $C^{1, \alpha}$ extension across the origin. 
From the results of \cite{DeTurckKazdan}, this is sufficient 
regularity to find a harmonic coordinate system around the origin. 
We view the equation as coupled to the equation for $g$ in harmonic 
coordinates:
\begin{align} 
\label{eee1}
\Delta Ric = Rm * Ric,\\
\label{eee2} 
\Delta g = Ric + Q(g, \partial g).
\end{align}
From (\ref{eee1}), as in \cite[Lemma 5.8]{BKN}, 
it is not hard to conclude that $Ric \in L^p$ for any $p < \infty$
(this is because from assumption we have a Sobolev constant bound,
and we also have an upper volume growth bound). 
Since $g$ is $C^{1,\alpha}$, 
from elliptic regularity, (\ref{eee2}) implies 
that $g \in W^{2, p}$, and therefore 
$Rm \in L^{p}$ for any $p$. Equation (\ref{eee1})
then implies $Ric \in W^{2,p}$, and (\ref{eee2}) 
gives $ g \in W^{3,p}$. Bootstrapping in this manner, 
we find that $g \in C^{\infty}$. 
\end{proof}
\section{Convergence}
\label{secconv}
In this section we complete the proofs of 
Theorem \ref{orbthm}, \ref{orbthm2}, and \ref{orbthm3}. 
We first describe the construction of the limit 
space, we will be brief since this step 
is quite standard (see for example \cite[Section 4]{Akutagawa}, 
\cite[Section5]{Anderson}, \cite[Section4]{Nakajima2}, 
\cite[Section 3]{Tian}). From the Sobolev constant bound (\ref{bound1})  
and lower volume bound (\ref{diamvol}), we obtain 
a lower growth estimate on volumes of geodesic 
balls. That is, there exists a constant $v > 0$ 
so that $Vol(B(x,s)) \geq v s^4$, for all $x \in M_j$
and $s \leq s_0$, for some $s_0 > 0$
\cite[Lemma 3.2]{Hebey}. 
For $R > 0$ large, let $M_{j,R} = M_j \cap B(p_j,R)$, 
and for $r > 0$ small, we take a maximally $r$-separated set
of $M_{j,R}$, that 
is, a collection of points $p_{i,j} \in M_{j,R}$ so that 
$B(p_{i,j}, r) \cap B(p_{i',j}, r) = \emptyset$ for $i \neq
i'$, and the collection $B(p_{i,j}, 2r)$ covers $M_{j,R}$. 
From the assumed bound (\ref{l2curv}) on the $L^2$-norm of curvature,
only a uniformly finite number of the balls $B(p_{i,j}, r)$ may satisfy
$ \int_{B(p_{i,j}, r)} |Rm_{g_j}|^2 dV_j \geq \epsilon_0$,
where $\epsilon_0$ is the constant from Theorem \ref{higherlocalregthm}.
By passing to a subsequence, we may assume that the 
number of these points is constant. Let us denote this collection 
of points by $S_j$, let $S_j(r)$ denoted the $r$-neighborhood of $S_j$,
and let $\Omega_j(r) = M_{j,R} \setminus S_j(r)$. 
From Theorem \ref{higherlocalregthm}, the curvature
and all covariant derivatives are uniformly bounded
on compact subsets of $\Omega_j(r)$. Furthermore,
the lower volume growth estimate implies an injectivity 
radius estimate (see \cite{CGT}), so we may apply a version of 
the Cheeger-Gromov convergence theorem (see \cite{Anderson}, \cite{Tian}) to find a 
subsequence such that 
$(\Omega_j(r), g_j)$ converges smoothly to 
$(\Omega_{\infty}(r), g_{\infty})$ as 
$j \rightarrow \infty$ on compact subsets.
That is, there exist diffeomorphisms $\Phi_{j,r} :
\Omega_{\infty}(r) \rightarrow \Omega_j(r)$ such 
that $\Phi_{j,r}^{\ast} g_j$ converges to $g_{\infty}$ in 
$C^{\infty}$ on compact subsets of $\Omega_{\infty}(r)$.  
By choosing a sequence $r_i \rightarrow 0$, and by 
taking diagonal subsequences, we obtain limit spaces with natural 
inclusions $\Omega_{\infty}(r_i) \subset \Omega_{\infty}(r_{i+1})$.
Letting $i \rightarrow \infty$, we obtain a limit 
space $(M_{\infty,R}, g_{\infty})$.
This is done for each $R$ large, and taking a sequence 
$R_i \rightarrow \infty$, we obtain a pointed limit 
space $(M_{\infty}, g_{\infty}, p_{\infty})$.

 We will now show how the main part of Theorem \ref{orbthm} follows 
assuming Theorem \ref{orbthm2}, and then we will complete the proof of 
Theorem \ref{orbthm2} below. 
In fact, we only require the volume growth 
estimate from Theorem \ref{orbthm2} to hold 
only for $r \leq r_0$, where $r_0$ is some fixed scale. 
That is, let us assume that 
\begin{align}
\label{volloc}
Vol(B_{g_i}(p,r)) \leq V r^4
\end{align}
 for all
$ p \in M_i$, and all $r \leq r_0$. 
The volume growth estimate (\ref{volloc})  
implies that we may add finitely many points to 
$M_{\infty}$ to obtain a complete length space;
this follows since $\#|S_j|$ is uniformly bounded
(see \cite[Section 5]{Anderson}, \cite[Section 3]{Tian}) for 
more details). For notational simplicity, we will continue to 
denote the completion by $M_{\infty}$.

 The estimate (\ref{volloc}), together with a global lower 
volume bound, imply a lower diameter
bound $\mbox{diam}(M_i, g_i) > \lambda >0$, which 
implies that $M_{\infty} \neq S$.
From (\ref{volloc}), it follows also that we have local 
volume convergence, and $(M_j, g_j, p_j)$ converges to
 $(M_{\infty}, g_{\infty}, p_{\infty})$ in the Gromov-Hausdorff
distance, moreover, the convergence is of length spaces.

 To analyze the singular points of $M_{\infty}$, for $p \in S$
we look at $B(p,\delta) \setminus \{p \}$ for $\delta$ small. 
The volume growth estimate (\ref{volloc}) implies the 
number of components of $B(p,\delta) \setminus \{p \}$
is finite (see \cite[Lemma 3.4]{Tian}). Restricting to each 
component, Theorem \ref{remsing} 
implies that the singularities are metric orbifold singularities, 
that is, the metric is locally a quotient of a 
{\em{smooth}} metric on each cone. 
Consequently, $M_{\infty}$ is a Riemannian multi-fold. 
Using what we have proved so far about limits 
(i.e., under the assumption (\ref{volloc})), we
next prove Theorem \ref{orbthm2}. 
\begin{proof}(of Theorem \ref{orbthm2}).
By Theorem \ref{higherlocalregthm}, if $g$ is critical, and 
\begin{align*}
\Vert Rm \Vert_{L^2(B(p,2))} = 
\left\{ \int_{B(p,2)} |Rm|^2 dV_g \right\}^{1/2} \leq \epsilon_0,
\end{align*}
then 
\begin{align*}
\underset{B(p, 1)}{sup}| Rm| \leq \frac{1}{4} C \epsilon_0. 
\end{align*}
By Bishop's volume comparison theorem, we must have
$Vol(B(p,1)) \leq A_0$, where $A_0$ depends only on the 
Sobolev constant.  

 We also note the following fact, for any metric, 
\begin{align*}
\lim_{r \rightarrow 0} Vol( B(p,r))r^{-4} = \omega_4,
\end{align*}
where $\omega_4$ is the volume ratio of the Euclidean metric
on $\mathbf{R}^4$. Clearly, $ A_0 \geq \omega_4$.

For any metric $(M,g)$, define the maximal volume ratio as
\begin{align} 
MV(g) = \underset{x \in M, r \in \mathbf{R}^+}{\max}
\frac{ Vol ( B(x,r))}{r^4}.
\end{align}
If the theorem is not true, then there exists 
a sequence of critical manifolds $(M_i,g_i)$, 
with $MV(g_i) \rightarrow \infty$, that is, 
there exist points $x_i \in M_i$, and $t_i \in \mathbf{R}^+$ so that 
\begin{align}
\label{contvol}
Vol (B(x_i,t_i))\cdot t_i^{-4} \rightarrow \infty, 
\end{align}
as $i \rightarrow \infty$. We choose a subsequence (which for 
simplicity we continue to denote by the index $i$) and radii 
$r_i$ so that 
\begin{align}
\label{2wchoice}
2 \cdot A_0 = Vol( B(x_i,r_i)) \cdot r_i^{-4} 
= \underset{r \leq r_i}{ \max} \ Vol( B(x_i,r)) \cdot r^{-4},
\end{align}
We furthermore assume that $x_i$ is chosen so that 
$r_i$ is minimal, that is, the smallest radius for which 
there exists some $p \in M_i$ such that 
$Vol(B_{g_i}(p, r)) \leq 2 A_0 r^4$, for all $r \leq r_i$.

First let us assume that $r_i$ has a subsequence converging to zero.
For this subsequence (which we continue to index by $i$), 
we consider the rescaled metric $\tilde{g}_i = r_i^{-2}g_i$, 
so that $B_{g_i}( x_i, r_i)  = B_{\tilde{g}_i}( x_i, 1)$. 
From the choice of $x_i$ and $r_i$, 
the metrics $\tilde{g}_i$ have bounded volume ratio,
in all balls of unit size.  

From the argument above, some subsequence converges
on compact subsets to a complete 
length space $(M_{\infty}, g_{\infty}, p_{\infty})$ 
with finitely many point singularities.  
The limit could either be compact or non-compact. 
In either case, the arguments above imply that the 
limit is a Riemannian multi-fold. 

 \begin{claim} The limit $(M_{\infty}, g_{\infty}, p_{\infty})$ contains 
at most $B_1$ disjoint bad annuli. 
\end{claim}
\begin{proof} 
We know that $(M_i, \tilde{g}_i, x_i)$ converges to $(M_{\infty}, g_{\infty}, 
p_{\infty})$
as pointed spaces. Assume that $(M_{\infty}, g_{\infty}, p_{\infty})$
contains $B_1 + 1$ disjoint bad annuli $A_l, 1 \leq l \leq B_1 + 1$.
Then there exists a radius $R$ so that $\cup A_l \subset B(p_{\infty}, R)$. 
Since the convergence is of pointed spaces, given any $\epsilon > 0$, 
there exist pointed, continuous $\epsilon$-almost isometries
$\Phi_{i, \epsilon} : B_{\tilde{g}_i}(x_i, 2R) \rightarrow 
B_{g_{\infty}}(p_{\infty}, 2R + \epsilon)$ for $i$
sufficiently large (see \cite{Burago}). 
For $\epsilon$ sufficiently small, it is easy to see that 
for each $l$, $\Phi_{i, \epsilon}^{-1}(A_l)$ will be $\epsilon$-close \
to a bad annulus 
in $(M_i, \tilde{g}_i, x_i)$. Applying the Mayer-Vietoris argument in 
\cite[Lemma 4.7]{TV} to this collection, we conclude that the number 
must be bounded by $B_1$, a contradiction. 
\end{proof}

If $M_{\infty}$ is non-compact, the remarks at the end of Section 
\ref{localregularity} shows that assumption (\ref{decay1}) is
satisfied. Also, from \cite[Lemma 6.1]{TV}, the Sobolev 
constant bound implies a lower volume growth bound 
(this is valid for orbifolds), so (\ref{cond4}) is satisfied. 
Theorem \ref{bigthm} then implies 
that $M_{\infty}$ has only finitely many ends, and that there 
exists a constant
$A_1 \geq 2 A_0$ so that 
\begin{align}
\label{poq}
Vol(B_{g_{\infty}}(p_{\infty}, r)) \leq 
A_1 r^4, \mbox{ for all } r > 0.  
\end{align}
  If $M_{\infty}$ is compact, then clearly the estimate 
(\ref{poq}) is valid for some 
constant $A_1 \geq 2 A_0$, since the limit has finite 
diameter and volume, and the estimate holds for $r \leq 1$.    

The inequality 
\begin{align}
\label{ballineq}
\int_{ B_{g_i}( x_i, r_i)} |Rm_i|^2 dV_i \geq \epsilon_0,
\end{align}
must hold; otherwise, as remarked above, we would have 
$Vol (B_{g_i}( x_i, r_i)) \leq A_0 r_i^4$, which 
violates  (\ref{2wchoice}). 

If the $r_i$ are bounded away from zero then 
there exists a radius $t$ so that 
\begin{align}
Vol(B_{g_i}(p,r)) \leq 2A_0 r^4, \mbox{ for all } r \leq t, p \in M_i.
\end{align}
We repeat the argument from the first case, but 
without any rescaling. 
Since the maximal volume ratio is bounded on small scales, 
we can extract an orbifold limit. The limit can 
either be compact or non-compact, but the inequality 
(\ref{poq}) will also be satisfied for some $A_1$, 
Following the same argument, we find a sequence of balls 
satisfying (\ref{ballineq}). 

We next return to the (sub)sequence $(M_i, g_i)$
and extract another subsequence so that 
\begin{align}
\label{3wchoice}
200 \cdot A_1 = Vol( B(x_i',r_{i}')) \cdot (r_i')^{-4} 
= \underset{r \leq r_i'}{ \max} \ Vol( B(x_i',r)) \cdot r^{-4}.
\end{align}
Again, we assume that $x_i'$ is chosen so that 
$r_i'$ is minimal, that is, the smallest radius for which 
there exists some $p \in M_i$ such that 
$Vol(B_{g_i}(p, r)) \leq 200 A_1 r^4$, for all $r \leq r_i$.
Clearly, $r_i < r_i'$.  

Arguing as above, if $r_i' \rightarrow 0$ 
as $ i  \rightarrow \infty$,  then we 
repeat the rescaled limit construction, but now 
with scaled metric $\tilde{g}_i = (r_i')^{-2} g_i$, 
and basepoint $x_i'$. We find a limiting 
orbifold $(M_{\infty}', g_{\infty}', p_{\infty}')$,
and a constant $A_2 \geq 200 A_1$ so that 
\begin{align*}
Vol(B_{g_{\infty}'}(p_{\infty}', r)) \leq A_2 r^4 \mbox{ for all } r > 0. 
\end{align*}
For the same reason as above, we must have
\begin{align*}
\int_{ B_{g_j}( x_j', r_j')} |Rm_j|^2 dV_j \geq \epsilon_0.
\end{align*}
If $r_i'$ is bounded below, we argue similarly, but without any rescaling. 

We claim that for $i$ sufficiently large, the 
balls $B(x_i,r_i)$ (from the first subsequence)
and $B(x_i', r_i')$ (from the second) must be disjoint 
because of the choice in (\ref{3wchoice}). 
To see this, if  $B(x_i,r_i) \cap B(x_i', r_i') \neq \emptyset$,
then $  B(x_i', r_i') \subset B(x_i, 3r_i')$.
Then (\ref{poq}) and (\ref{3wchoice}) yield
\begin{align*}
200 A_1 (r_i')^4 = Vol ( B ( x_i', r_i'))
\leq Vol (  B ( x_i, 3 r_i')) \leq 2 A_1 ( 3 r_i')^4
 = 162 A_1 (r_i')^4,
\end{align*}
which is a contradiction
(note the last inequality is true for $i$ sufficiently large 
since (\ref{poq}) holds for the limit).

We repeat the above procedure. The process must terminate in 
finitely many steps from the bound 
$ \Vert Rm_i \Vert_{L^2} < \Lambda$. This contradicts 
(\ref{contvol}), which finishes the proof.

\end{proof}

The convergence statement in Theorem \ref{orbthm}
now follows from Theorem \ref{orbthm2}, since (\ref{volloc})
is satisfied. Statement (3) follows from 
Theorem \ref{decayale}, since the multi-fold is
the union of irreducible orbifolds. Note also that the 
volume bound in Theorem \ref{orbthm2} gives a uniform bound for the 
the number of irreducible pieces, and the number of ends of the limit
multi-fold. 

  To finish the proof of Theorem \ref{orbthm} we need
to verify statements (2), (4), and (5). 
The next proposition gives a direct argument to bound the 
number of components of $B(p,\delta) \setminus \{p \}$ for $\delta$ small
in terms of the Sobolev constant and first Betti number. 

\begin{proposition}
\label{finitely1}
For $p \in M_{\infty}$, and $\delta$ sufficiently small,  
the number of components of $B(p, \delta) \setminus \{p \}$ 
can be estimated in terms of the first Betti number and
the Sobolev constant 
(defined as in (\ref{mainSob} or \ref{mainSob2})). 
If $b_1(M_i) = 0$, then $p$ is irreducible. 
Furthermore, $M_{\infty}$ does not split off any compact orbifold. 

 If the weaker Sobolev inequality (\ref{weakSob}) is 
assumed, then the number of components of $B(p, \delta) \setminus \{p \}$ 
can be still be estimated in terms of the Sobolev constant and the 
first Betti number (but in this case it is possible that 
if $b_1(M_i) =0$, a singular point could be reducible, and 
it is also possible that $M_{\infty}$ could split off a compact orbifold).  
\end{proposition}
\begin{proof}
Let $p$ be a non-irreducible 
singular point. 
We have shown around $p$, $M_{\infty}$ is a finite
union of orbifold cones, with the basepoints 
identified. For each orbifold, since the convergence is smooth 
away from the singular points, we look before the limit,
and this gives us a portion of a cone 
on $S^3 / \Gamma$ in the original 
manifold, very small, which we call
$N_i \subset M_i$ and $N_i 
= (a_i, 2a_i) \times (S^3 / \Gamma),$
which is close, in any $C^{m}$-topology, 
to an annulus $A(a_i, 2a_i)$
in a cone on a spherical space form 
$\mathcal{C} (S^3 / \Gamma)$, and 
where $a_i \rightarrow 0$ as $i \rightarrow \infty$. 

If $ \{ a_i \} \times S^3 / \Gamma $
bounds a region in $M_i$, equivalently, if $N_{i}$ 
separates $M_i$ into two components, then 
this decomposes $M_i$ into a disjoint union $A_i \cup N_{i} \cup B_i$.
Since the point $p$ is non-irreducible and the convergence is smooth 
away from the singular points, we must have $ Vol(A_i)$ and $Vol(B_i)$  
uniformly bounded away from zero. Without loss of generality, 
assume $Vol(A_i) \leq Vol(B_i)$.  

  We take a function $f_i$ which is $1$ on $A_i$, $0$ on $B_i$, 
since the neck $N_{i}$ is $C^{m}$-close to the annulus
$A(a_i,2a_i)$ in a flat cone, we may take $|\nabla f| = 1/a_i$ on 
the neck portion $N_{i}$. 
We compute 
\begin{align}
\Vert f_i \Vert_{L^{4}}  
= \left\{ \int_{A_i} 1 dV_{g_i} + \int_{N_i} f_i dV_{g_i}
\right\}^{1/4}
\sim Vol(A_i)^{1/4}.
\end{align} 
Next, 
\begin{align}
\Vert \nabla f_i \Vert_{L^2}^2
 = \int_{N_i} \frac{1}{a_i} dV_{g_i} 
\sim \frac{1}{a_i} C ( (2a_i)^4 - (a_i)^4) = C a_i^3,
\end{align}
since $N_i$ is $C^{m}$-close to $A(a_i,2a_i)$. 
Using the Sobolev inequality (\ref{mainSob}), we obtain
\begin{align}
  Vol(A_i)^{1/4}  \leq C_S 
  C' a_i^{3/2}  + Vol(M_i)^{-1/4} Vol(A_i)^{1/2},
\end{align} 
Rearranging terms, 
\begin{align}
\label{rearr}
 Vol(A_i)^{1/4}( 1  -  Vol(M_i)^{-1/4} Vol(A_i)^{1/4} )  \leq C_S 
  C' a_i^{3/2}. 
\end{align} 
We have $Vol(M_i) \geq 2 Vol(A_i)$, therefore 
\begin{align}
 Vol(A_i)^{1/4}( 1  -  2^{-1/4})   \leq C_S 
  C' a_i^{3/2}. 
\end{align}
Since $ Vol(A_i)$ is uniformly bounded 
away from zero, this is a contradiction for $i$ large. 
Therefore none of the necks $N_i$ around a non-irreducible 
singular point can bound regions in $M_i$. Using the 
intersection pairing, any of these 
embedded space forms will give a generator 
of $b_1(M_i)$. At most two of these may give
rise to the same generator, so from the 
assumed bound on $b_1(M_i)$, there may only 
be finitely many, and if $b_1(M_i) =0$, the
singular point $p$ must be irreducible. 

  Note that in case of the Sobolev inequality (\ref{mainSob2}), a
similar argument works, and a similar argument 
shows that $M_{\infty}$ does not split 
off any compact orbifold.

In the case (\ref{weakSob}) is satisfied, 
let $p$ be non-irreducible 
singular point. Again, we have shown around 
$p$, $M_{\infty}$ is a finite union of orbifold cones, 
with the basepoints identified.
For each orbifold group $\Gamma_j$ at $p$, since the convergence is smooth 
away from the singular points, we look before the limit,
and this gives us a portions of cones 
on $S^3 / \Gamma_{j}$ in the original 
manifold, $N_{i,j} \subset M_i$, very small, $N_{i,j} 
= (a_i, 2a_i) \times (S^3 / \Gamma_j),$
which is close, in any $C^{m}$-topology to an 
annulus $A(a_i, 2a_i)$
in a cone on a spherical space form 
$\mathcal{C} (S^3 / \Gamma_j)$, and 
where $a_i \rightarrow 0$ as $i \rightarrow \infty$. 

 Take any collection of $Q > 16 C_S^4$ 
irreducible orbifolds at $p$. Then we claim 
at least one of the necks $N_{i,j}$ cannot bound
a region in $M_i$, i.e., $N_{i,j}$ cannot separate 
$M_i$ into 2 components. 
 If all of the $N_{i,j}$ bound, then 
this decomposes $M_i$ into a disjoint union 
$A_i \cup ( \cup_{j} N_{i,j}) \cup B_i$,
where $A_i$ is taken to be on the side of the neck
where convergence is smooth, $B_i$ is the 
rest of $M_i$. Since we have a finite 
collection, and convergence is smooth on $A_i$, 
so $ Vol(A_i)$ is uniformly bounded away from zero. 
Now $A_i$ is the union of $Q$ regions, therefore, 
one of the regions, which we call $R_{i,j}$, 
must satisfy $Vol(R_{i,j}) < \frac{1}{Q} Vol(A_i)$. 

  We take a function $f_i$ which is $1$ on the region 
$R_{i,j}$, 
since the neck $N_{i,j}$ bounding $R_{i,j}$ 
is $C^{\infty}$-close to the annulus
$A(a_i,2a_i)$ in a flat cone, we may take $|\nabla f| = 1/a_i$ on 
the neck portion $N_i$, with $f_i =0 $ otherwise.

 As in (\ref{rearr}) above, but using the 
Sobolev inequality (\ref{weakSob}), we obtain 
\begin{align}
 Vol(R_{i,j})^{1/4}( 1  -  C_S Vol(M_i)^{-1/4} Vol(R_{i,j})^{1/4} )  \leq C_S 
  C' a_i^{3/2}. 
\end{align} 
We have $Vol(M_i) \geq Q Vol(R_{i,j})$, therefore 
\begin{align}
 Vol(R_{i,j})^{1/4}( 1  - C_s Q^{-1/4})   \leq C_S 
  C' a_i^{3/2}, 
\end{align}
from the choice of $Q$, we obtain
\begin{align}
\frac{1}{2}  Vol(R_{i,j})^{1/4} \leq C_S 
  C' a_i^{3/2}.
\end{align}
Since $ Vol(R_{i,j})$ is uniformly bounded 
away from zero, this is a contradiction for $i$ large. 
Therefore, for any collection of $Q > 16 C_S^4$ 
irreducible orbifolds at $p$, one 
of the neck $N_{i,j}$ cannot bound 
a regions in $M_i$. Using the 
intersection pairing, the corresponding 
embedded space form $S^3 / \Gamma_{i,j}$ 
will give a generator 
of $b_1(M_i)$. 

If there are $k*Q$ orbifolds at $p$, then we 
find $k$ generators  $b_1(M_i)$.
At most $2$ of these may give
rise to the same generator, so from the 
assumed bound on $b_1(M_i)$, there may only 
be finitely many. 
\end{proof}

 We remark that we may characterize the singular set 
as follows: with $\epsilon_0$ as in 
Theorem \ref{higherlocalregthm}, we have 
\begin{align}
\begin{split}
S = \{& x \in M_{\infty} | \liminf_{j \rightarrow \infty}
\int_{B(x_j, r)} |Rm_{g_j}|^2 dvol_{g_j} \geq \epsilon_0 \\
& \mbox{for any sequence } \{x_j \} \mbox{ with }  \lim_{j \rightarrow
\infty} x_j = x, \mbox{ and all } r > 0 \}.
\end{split}
\end{align} 

We next give a description of the convergence 
at the singular points, by rescaling the sequence 
at a singular point $x \in S$. 
Several bubbles may arise in the degeneration, so we have 
to rescale properly, and possibly at several different 
scales. This was done in \cite{Bando} for the 
Einstein case, and with a few minor changes, 
the proof works in our case. We outline 
the details here. 
 For $0 < r_1 < r_2$, we let $D(r_1,r_2)$ denote 
$B(p, r_2) \setminus B(p, r_1)$.
Given a singular point $x \in S$, we take a 
sequence $x_i \in (M, g_i)$ such that 
$\lim_{ i \rightarrow \infty} x_i = x$
and $B(x_i, \delta)$ converges to 
$B(x, \delta)$ for all $\delta > 0$. 
 We choose a radius $r_{\infty}$  sufficiently 
small and the sequence $x_i$ to satisfy
\begin{align}
\underset{B( x_i, r_{\infty})}{ \mbox{sup}} |Rm_{g_i}|^2 = 
|R_{g_i}|^2( x_i) \rightarrow \infty \mbox{ as } j \rightarrow
\infty,  
\end{align}
and 
\begin{align} \int_{B(x,r_{\infty})} |R_{g_{\infty}}|^2 dV_{g_{\infty}}
  \leq \epsilon_0 /2.
\end{align}

We next choose $r_0(j)$ so  that 
\begin{align}
 \int_{D(r_0, r_{\infty})} |R_{g_j}|^2 dV_{g_j} = \epsilon_0, 
\end{align}
where $\epsilon_0$ is as in Theorem \ref{higherlocalregthm}, 
and again $D(r_o, r_{\infty}) = B(x_i, r_{\infty}) \setminus 
B(x_i, r_0)$. 
An important note, which differs from the 
Einstein case, the annulus $D(r_0, r_{\infty})$ 
may have several components. 

 Since the curvature is concentrating at $x$, 
$r_o(j) \rightarrow 0$ as $j \rightarrow \infty$. 
From Theorem \ref{orbthm}, the rescaled sequence 
$(M, r_o(j)^{-2} g_i, x_i)$ has a subsequence 
which converges to a complete, non-compact 
multi-fold with finitely many singular points, 
which we denote by $M_{i_1}, 1 \leq i_1 \leq 
\# \{S\}$. 
Since 
\begin{align}
\int_{D(1, \infty)} |Rm|^2 dV_g \leq \epsilon_0, 
\end{align}
there are no singular points outside 
of $B(x,1)$. 

 On the noncompact ends, from Theorem \ref{decayale},
the metric is ALE  of order $\tau$ for any $\tau < 2$.
As in \cite[Proposition 4]{Bando}, we conclude 
that the neck region (for large $i$) will be 
arbitrarily close to a portion of a flat cone 
$\mathbf{R}^4 / \Gamma$, possibly several cones 
if $M_{i_1}$ has several ends. 
So the convergence at a singular point $x_{i_1}$ is 
that the ALE multi-fold $M_{i_1}$ is 
bubbling off, or scaled down to a point 
in the limit.

 To further analyze the degeneration at the 
singular points, we look at the multi-fold $M_{i_1}$
with singular set $S_{i_1}$. If $S_{i_1}$ is
empty, then we stop. We do the 
same process as above for each singular 
point of $M_{i_1}$, and obtain ALE multi-folds 
$M_{i_1, i_2}, 1 \leq i_2 \leq \# \{S_{i_1} \}$. 
If $M_{i_1, i_2}$ has singularities, then 
we repeat the procedure. This process 
must terminate in finite steps, since in this 
construction, each singularity takes at least
$\epsilon_0$ of curvature. As pointed out in 
\cite{Bando2}, there could be some 
overlap if any singular point lies 
on the boundary of $B(1)$ at some stage 
in the above construction. But there can 
only be finitely many, and then there 
must also be a singular point in the 
interior of $B(1)$, so we still take 
away at least $\epsilon_0$ of curvature 
at each step. 

 In the K\"ahler case, one can use 
the methods of \cite{Li-Tam}  to show that 
boundary of sufficiently small balls around the 
singular points of $M_{\infty}$ are connected.
If a singular point $p \in M_{\infty}$ is 
non-irreducible, then using the above 
bubbling analysis, at some step one must find 
an irreducible K\"ahler ALE orbifold with more than one 
end. From the remark following Theorem \ref{decayale}, 
this is not possible, therefore in the 
K\"ahler case (c), the limit is
irreducible. This completes the proof of Theorem \ref{orbthm}. 

\subsection{Local Sobolev inequality}

 As we have noted throughout the paper, many of our
results hold with a weaker assumption on the 
Sobolev constant. We have the following notion 
of local Sobolev constant. 
We define $C_{S}(r)$ to be the best 
constant such that  
\begin{align}
\label{mainSob3}
\Vert f \Vert_{L^{4}} \leq C_{S}(r)   
\Vert \nabla f \Vert_{L^2},
\end{align}
for all $f \in C^{0,1}$ with compact support 
in $B(p, r)$, and for all $p \in M$. 

 The following is the analogue of Theorem \ref{orbthm3} 
with a local Sobolev constant bound (the proof is identical):
\begin{theorem} 
\label{orbthmsob}
Let $(M_i, g_i)$ be a sequence of critical metrics  $g_i$ on 
smooth, four-dimensional manifolds $M_i$ satisfying
\begin{align} 
\label{sbound1}
&C_{S}(r_0) \leq C_1 \ \ \ (\mbox{for some fixed } r_0 > 0), \\
\label{sl2curv}
&\int_{M_i} |Rm_{g_i}|^2 dV_{g_i} \leq \Lambda,\\
\label{sdiamvol}
&  Vol(g_i) = 1,\\
\label{sbetti1}
& b_1(M_i) < B_1, 
\end{align}
where $C_1, \Lambda, \lambda$  are constants,
and $b_1(M_i)$ denotes the first Betti number,
Then there exists a subsequence $\{ j \} \subset \{ i \}$, 
a compact, connected, critical Riemannian multi-fold $(M_{\infty}, g_{\infty})$, 
and a finite singular set $S \subset M_{\infty}$ such that 
$(M_j, g_j)$ converges to $(M_{\infty}, g_{\infty})$.
\end{theorem}
\section{Further Remarks}
We conclude by listing here some interesting problems.
\vspace{2mm}
\\
1) We considered above the case of constant scalar curvature 
K\"ahler metrics. We conjecture that these results 
extend to the more general extremal K\"ahler case in 
dimension four (\cite{Calabi1}, \cite{Calabi2}). 
\vspace{2mm}
\\ 
2) It is an interesting problem to generalize 
our results to higher dimensions. We conjecture that 
the following is true for the higher dimensional 
extremal K\"ahler case. 
Assuming fixed K\"ahler class, first and second Chern 
classes, the limit space has at most a codimension $4$ singular set, 
and the singular set is a holomorphic subvariety. 
Even in the case of Bach-flat or harmonic
curvature in higher dimensions,  
under the bound $\Vert Rm \Vert_{L^2} < \Lambda$, the
limit space should have a most a codimension $4$ singular set, 
with top strata modeled on orbifold singularities. 
This was proved for Einstein metrics in \cite{CC1}, 
\cite{CC2}, \cite{CC3}, \cite{CCT}. 
\vspace{2mm}
\\ 
3) It would be very interesting to remove the 
Sobolev constant assumption and understand the collapsing 
case. 
\vspace{2mm}
\\ 
4) In the general Bach-flat case in dimension four,  
one should be able to show that the orbifold singularities 
are {\em{smooth}} metric singularities, and that in the ALE case, 
one can obtain a positive order of decay. 
\bibliography{Moduli_references}
\end{document}